\documentstyle[amssymb,amsfonts]{amsart}

\renewcommand\P{{\mathbf{P}}}
\newcommand\Q{{\mathbf{Q}}}
\newcommand\R{{\mathbf{R}}}
\newcommand\Z{{\mathbf{Z}}}
\newcommand\rank{{\operatorname{rank}}}
\newcommand\Image{{\operatorname{Image}}}

\newcommand\size{{\operatorname{size}}}
\newcommand\dist{{\operatorname{dist}}}
\newcommand\mes{{\operatorname{mes}}}
\renewcommand\mod{{\ \operatorname{mod}\ }}

\newcommand\eps{\varepsilon}
\def\endprf{\hfill  {\vrule height6pt width6pt depth0pt}\medskip}
\newenvironment{proof}{\noindent {\bf Proof} }{\endprf\par}

\theoremstyle{plain}
  \newtheorem{theorem}[subsection]{Theorem}

  \newtheorem{lemma}[subsection]{Lemma}
  \newtheorem{corollary}[subsection]{Corollary}
   
\theoremstyle{remark}
  \newtheorem{remark}[subsection]{Remark}
  
\theoremstyle{definition}
  \newtheorem{definition}[subsection]{Definition}
\include{psfig}
\begin{document}
\title[John-type theorems for GAPs and sumsets]{ John-type theorems for generalized arithmetic progressions and iterated sumsets }

\author{Terence Tao}
\address{Department of Mathematics, UCLA, Los Angeles CA 90095-1555}
\email{tao@@math.ucla.edu}
\thanks{T. Tao is supported by a grant from the Macarthur Foundation.}
\author{  Van  Vu}
\address{Department of Mathematics, Rutgers, Piscataway, NJ 08854}
\email{vanvu@@math.rutgers  .edu}
\thanks{V.   Vu  is supported by  NSF Career Grant 0635606.}


\begin{abstract} A classical theorem of Fritz John allows one to describe a convex body, up to constants, as an ellipsoid.  In this article we establish similar descriptions for generalized (i.e. multidimensional) arithmetic progressions in terms of proper (i.e. collision-free) generalized arithmetic progressions, in both torsion-free and torsion settings.  We also obtain a similar characterization of iterated sumsets in arbitrary abelian groups in terms of progressions, thus strengthening and extending recent results of Szemer\'edi and Vu.
\end{abstract}
\maketitle

\section{Introduction}

Define a \emph{convex body} to be a compact convex subset of a Euclidean space $\R^d$ with non-empty interior\footnote{This differs slightly from the notation in \cite{TVbook}, in which convex bodies were assumed to be open and bounded rather than compact.  This change is convenient for some minor technical reasons, but does not significantly affect any of the results given here.}. We say that a convex body $B$ is \emph{symmetric} if $-1 \cdot B = B$, where $\lambda \cdot B := \{\lambda x: x \in B\}$ denotes the dilate of $B$.
A classical theorem of John \cite{john} characterizes such bodies up to constants:

\begin{theorem}[John's theorem, symmetric case]\label{john-thm}  Let $B$ be a symmetric convex body in $\R^d$.  Then there exists a (closed) ellipsoid $E$ centered at the origin such that
$$ E \subseteq B \subseteq \sqrt{d} \cdot E.$$
\end{theorem}

The constant $\sqrt{d}$ here is sharp, as can be seen by considering the case when $B$ is a box or cube.  There is an analogous theorem for asymmetric convex bodies, but we will consider mainly symmetric situations here.

The purpose of this paper is to investigate discrete analogues of John's theorem, when $\R^d$ is replaced by a lattice $\Gamma$, a progression, or an iterated sumset.  or more generally an arbitrary additive group $G$.  For this, we shall need to replace the concept of an ellipsoid by the notion of a (proper) \emph{generalized arithmetic progression} (GAP), which we now pause to define.  Again we
restrict attention to the symmetric case.

\begin{definition}[Sumset notation] An \emph{additive group} is an abelian group $G = (G,+)$.
If $A,B$ are sets in an additive group $G$, we use $A+B := \{ a+b: a \in A, b \in B \}$ for the sumset, $A-B := \{a-b: a \in A, b \in B \}$ for the difference set, and $|A|$ for the cardinality of $A$.  For $k \geq 1$, we use $kA = A+ \ldots + A$ to denote the iterated sumset of $k$ copies of $A$.  If $n$ is an integer, we use $n \cdot A := \{ na: a \in A \}$ to denote the dilate of $a$ by $n$, where $na$ is the sum of $n$ copies of $a$ (with the conventions $0a =0$ and $(-n)a = -(na)$).  We caution that $k \cdot A \neq kA$ in general, although we do have $k \cdot A \subseteq kA$.  If $I$ is a set of integers, we write $I \cdot A := \{ n a: n \in I, a \in A \}$.
If $a,b$ are reals, we use $[a,b]_\Z := \{ n \in \Z: a \leq n \leq b \}$ to denote the discrete interval and $[a,b]_\R := \{ x \in \R: a \leq x \leq b \}$ to denote the continuous interval.  Similarly define $[a,b)_\Z$, $(a,b)_\R$, etc.
\end{definition}

\begin{definition}[GAPs]  Let $G$ be an additive group.  A \emph{symmetric generalized arithmetic progression} in $G$, or \emph{symmetric GAP} for short, is a triplet $\P = (N, v, d)$, where the \emph{rank} $\rank(\P) = d$ is a non-negative integer, the \emph{dimensions} $N = (N_1,\ldots,N_d)$ are a $d$-tuple of positive reals, and the \emph{steps} $v = (v_1,\ldots,v_d)$ are a $d$-tuple of elements of $G$.
\begin{itemize}
\item We define the \emph{image} $\Image(\P) \subset G$ of $\P$ to be the set
\begin{align*}
\Image(\P) &:= [-N_1,N_1]_\Z \cdot v_1 + \ldots + [-N_d,N_d]_\Z \cdot v_d \\
&=
\{ \sum_{i=1}^d n_i v_i: n_i \in [-N_i,N_i]_\Z \forall i =1,\ldots,d\}.
\end{align*}
\item For any $t > 0$, we define the \emph{dilate} $\P_t$ of $\P$ to be the GAP $\P_t := (t N, v, d)$ formed by dilating all the dimensions by $t$.
\item We say that $\P$ is \emph{proper} if all the elements $n_1 v_1 + \ldots + n_d v_d$ for $n_i \in [-N_i,N_i]_\Z$ are distinct.  More generally, we say that $\P$ is \emph{$t$-proper} for some $t > 0$ if $\P_t$ is proper, and \emph{infinitely proper} if it is $t$-proper for all $t$ (i.e. the elements $v_1,\ldots,v_d$ are algebraically independent).
\item We define the \emph{size} of $\P$ to be $\size(\P) := |\Image(\P)|$.   Observe that if $\P$ is proper if and only if
$\size(\P) = \prod_{i=1}^d (2 \lfloor N_i \rfloor + 1)$.
\end{itemize}
\end{definition}

\begin{remark} For technical reasons, we need to allow the dimensions $N_i$ to be real rather than integer.
Because of this, it is not always the case that $\Image(\P_t) + \Image(\P_{t'})$ is equal to $\Image(\P_{t+t'})$, although this is true if the components of $tN$ and $t'N$ have fractional parts in
$[0,1/2)_\R$.  Instead, we only have the inclusion $\Image(\P_t) + \Image(\P_{t'}) \subseteq \Image(\P_{t+t'})$ in general.
On the other hand, by replacing each of the dimensions with their greatest integer part we can always assume
that the dimensions are integer without affecting the image, rank or properness of $\P$ (although the image and properness of the dilates $\P_t$, will be affected).
\end{remark}

\begin{remark}
In most treatments, the progression $\P$ is identified with its image $\Image(\P)$.  However we shall avoid doing this here because many important features of the progression (such as the rank, or the dilates $\Image(\P_t)$) are not completely determined by the image alone.  In particular, if $\P$ and $\Q$ are symmetric GAPs, an inclusion $\Image(\P) \subseteq \Image(\Q)$ does not necessarily entail an inclusion $\Image(\P_t) \subseteq \Image(\Q_t)$ even when $t=2$ (unless the fractional parts of the $N_i$ all lie in $[0,1/2)_\R$).
\end{remark}

A \emph{lattice} is a discrete additive subgroup of a Euclidean
space. To start, we present an analogue of John's theorem on
lattices.

\begin{theorem}[Discrete John's theorem]\label{djt}  Let $B$ be a convex symmetric body in $\R^d$, and let $\Gamma$ be a lattice in $\R^d$.  Then there exists a symmetric, infinitely proper GAP $\P$ in $\Gamma$ with rank $\rank(\P) \leq d$ such that we have the inclusions
\begin{equation}\label{inclusions}
(O(d)^{-3d/2} \cdot B) \cap \Gamma \subseteq \Image(\P) \subseteq B \cap \Gamma \subseteq \Image(\P_{O(d)^{3d/2}})
\end{equation}
and more generally
\begin{equation}\label{inclusions-2}
 (O(d)^{-3d/2}t \cdot B) \cap \Gamma \subseteq \Image(\P_t) \subseteq (t \cdot B) \cap \Gamma
\end{equation}
for any $t > 0$.
Furthermore, we have the size bounds
\begin{equation}\label{dgam}
 O(d)^{-7d/2} |B \cap \Gamma| \leq \size(\P) \leq |B \cap \Gamma|.
\end{equation}
\end{theorem}

As usual $O(X)$ denotes a quantity bounded above by $CX$ for some absolute constant $C$; thus for instance $O(d)^{-d}$ denotes a quantity bounded from \emph{below} by $(Cd)^{-d}$ for some absolute constant $C$.

This theorem was essentially already established in \cite[Lemma
3.36]{TVbook}. In an earlier paper \cite{BV}, it was proved (see
Theorem 3 of \cite{BV}) that if $\Gamma$ is full dimensional then
one can find a GAP $\P$ such that $B $ is contained in the convex
hull $conv( \P)$ of $\P$ and $Vol (conv (\P)) \le C_d (Vol (B))$,
for some constant $C_d$ depending only on $d$.

For the convenience of the reader (and because this theorem will be
used to prove our other results) we supply a proof of Theorem
\ref{djt} in Section \ref{djt-sec}, taking the opportunity to
strengthen the bounds slightly and correct some misprints.

\subsection{The torsion-free case}

Next, we consider progressions in torsion-free additive groups $G$ (thus $nx \neq \{0\}$ for all $x \in G \backslash \{0\}$ and $n \in \Z \backslash 0$); for instance any lattice is torsion-free.  The natural question here is whether one can compare
a non-proper progression $\P$ with a proper or a $t$-proper progression.  In this regard, the following
results are known:

\begin{theorem}[Progressions contain proper progressions]\label{has-proper}\cite[Theorem 3.38]{TVbook}  Let $\P$ be a symmetric GAP of rank at most $d$ in an additive group $G$ (not necessarily torsion free).  Then there exists a proper symmetric GAP $\Q$ with $\rank(\Q) \leq \rank(\P)$, $\Image(\Q) \subseteq \Image(\P)$, and
$$ O(d)^{-5d} \size(\P) \leq \size(\Q) \leq \size(\P).$$
\end{theorem}

\begin{theorem}[Progressions contained in proper progressions]\label{in-proper}\cite[Theorem 3.40]{TVbook}, \cite[Theorem 2.1]{green} Let $\P$ be a symmetric GAP in a torsion-free group $G$, and let $t \geq 1$.  Then there exists a $t$-proper symmetric GAP $\Q$ with $\rank(\Q) \leq \rank(\P)$,
$\Image(\P) \subseteq \Image(\Q)$, and
$$ \size(\P) \leq \size(\Q) \leq (2t)^d d^{6d^2} \size(\P).$$
Furthermore, if $\P$ is not proper, we may take $\rank(\Q) \leq \rank(\P)-1$.
\end{theorem}

We remark that similar results were also obtained earlier by Bilu \cite{bilu} and Chang \cite{chang}.
The precise bound $(2t)^d d^{6d^2}$ is established in \cite{green}; the argument in \cite{TVbook} only considers the case $t=1$ and gives the weaker bound of $d^{O(d^3)}$.  We will be able to improve Theorem \ref{in-proper} in
Corollary \ref{prog-john2} (and Corollary \ref{prog-john-coset2}) below.

While these two results are already useful, they are not quite analogous to Theorem \ref{john-thm} or Theorem \ref{djt} because they only provide \emph{one-sided} containments; the original GAP $\P$ either contains or is contained in a proper GAP $\Q$, but the two proper GAPs on either side of $\P$ are not related to each other by a dilation.

Ideally we would like to take any GAP $\P$ of rank $d$ and a
parameter $t > 0$ and obtain inclusions $\Image(\Q_\eps) \subseteq
\Image(\P) \subseteq \Image(\Q)$ for some $t$-proper GAP $\Q$ and
some $\eps > 0$ depending in a reasonable manner on $d$ and $t$.
However this is not feasible once $t \geq 2$, even when $G$ is
just the integers $\Z$.  To see this, let $N$ be a large even
integer and consider the rank $2$ progression $\P := ((N/2,N),
(1,N), 2)$. The image of this progression fills out most of the
discrete interval $(-N^2 + N/2, N^2 - N/2)_\Z$, but misses the
half-integer multiples of $N$.  Suppose that we could find a
$2$-proper progression $\Q = (M, v, d)$ whose image contained
$\Image(\P)$.  Thus we have a map
$$f: \Image(\P) \to \prod_{i=1}^d [-M_i,M_i]_\Z \subset \Z^d $$
such that $x = f(x) \cdot v$ for all $x \in \Image(\P)$.  Since $\Q$ is $2$-proper, we easily see that $f$ is not only uniquely defined, but is also a \emph{Freiman homomorphism}, thus
$$ f(x_1) + f(x_2) = f(x_3) + f(x_4) \hbox{ whenever } x_1,x_2,x_3,x_4 \in \Image(\P) \hbox{ and } x_1 + x_2 = x_3 + x_4.$$
From this and the explicit form of $\Image(\P)$ it is not hard to show that $f$ must now be linear,
thus $f(x) = x f(1)$.  In particular, the discrete box $\prod_{i=1}^d [-M_i,M_i]_\Z^d$ must contain the one-dimensional progression $[-N^2+N/2+1,N^2-N/2-1]_\Z \cdot f(1)$.  From this we see that $\Image(\Q_\eps)$ contains
$[\eps (-N^2+N/2+1),\eps (N^2-N/2-1)]_\Z$, and so for fixed $\eps$ and large $N$ we cannot have $\Image(\Q_\eps) \subseteq \Image(\P)$.

It is likely that similar examples can also be constructed to cover the case $t=1$; the GAP $\P$ would thus be almost proper (with only a small number of collisions), and its image would cover most of the image of a genuinely proper GAP $\Q$, but there would be enough ``holes'' in the image that $\P$ could not accomodate a smaller version $\Q_\eps$ of $\Q$ for fixed $\eps$.

Nevertheless, we can get around this problem by weakening the conclusion in a number of ways.  One way is to no longer require the outer progression to be proper, but instead demand only that the inner progression is proper.  This leads to our first main result:

\begin{theorem}[John's theorem for GAPs]\label{prog-john} Let $\P$ be a symmetric GAP of rank $d \geq 0$ in a torsion-free group $G$, and let $t \geq 1$.  Then there exists a $t$-proper symmetric GAP $\Q$ of rank at most $d$,
such that we have the inclusions
$$ \Image(\P_{t'}) \subseteq \Image( \Q_{O(d)^{3d/2} t t'} )$$
for all $t' > 0$ and
$$ \Image(\Q_{tt'}) \subseteq \Image(\P_{t'})$$
for all $t' \geq 1$.
In particular, we have
$$ \Image(\Q) \subseteq \Image(\P) \subseteq \Image( \Q_{O(d)^{3d/2} t} ).$$
Furthermore we have the size bound
$$ t^{-d} 2^{-d^2-O(d \log d)} \size(\P) \leq \size(\Q) \leq \size(\P).$$
Finally, if $\P$ is not $1/2$-proper, then $\Q$ can be chosen to have rank at most $d-1$.
\end{theorem}

We prove this theorem (a fairly simple consequence of Theorem \ref{djt} and a rank reduction
argument) in Section \ref{pjt-sec}.  Note that it implies an analogue of Theorem \ref{has-proper} but with significantly weaker constants.

As a corollary we can obtain the following improvement of Theorem \ref{in-proper}.

\begin{corollary}[John's theorem for GAPs, again]\label{prog-john2}
Let $\P$ be a symmetric GAP of rank $d \geq 0$ in a torsion-free group $G$, and let $t \geq 1$.  Then there exists a $t$-proper symmetric GAP $\Q$ of rank at most $d$,
such that we have the inclusions
$$ \Image(\P_{t'}) \subseteq \Image( \Q_{t'} )$$
for all $t' > 0$ and
$$ \Image(\Q_{tt'}) \subseteq \Image(\P_{O(d)^{3d/2} t t'})$$
for all $t' \geq 1$.  In particular, we have
$$ \Image(\P) \subseteq \Image(\Q) \subseteq \Image( \P_{O(d)^{3d/2} t} ).$$
We also have the size bound
$$ \size(\P) \leq \size(\Q) \leq t^d O(d)^{3d^2/2} \size(\P).$$
Finally, if $\P_{d^{3d/2} t}$ is not proper, then $\Q$ can be chosen to have rank at most $d-1$.
\end{corollary}

\begin{proof} Apply Theorem \ref{prog-john} with $\P$ replaced by $\P_{O(d)^{3d/2} t}$ for a sufficiently large choice of $O(d)$.
\end{proof}

\subsection{The torsion case}

Now we turn to the case where $G$ is allowed to contain torsion (in particular, $G$ could be a finite group); equivalently, $G$ contains non-trivial finite subgroups.  Here, it is no longer reasonable to work with $t$-proper GAPs for any $t \geq 2$.  For instance, if $G$ is a non-trivial finite group, then (by the classification of such groups) $G$ is the image of a proper GAP, but cannot be the image of a $t$-proper GAP for any $t \geq 2$.  Instead, as first observed by Green and Ruzsa \cite{GR}, one should replace GAPs by the more general notion of a \emph{coset progression}:

\begin{definition}[Coset progressions]  Let $G$ be an additive group.  A \emph{symmetric coset progression} in $G$ is a quadruplet $\P = (N, v, d, H)$, where the \emph{rank} $\rank(\P) = d$ is a non-negative integer, the \emph{dimensions} $N = (N_1,\ldots,N_d)$ are a $d$-tuple of positive reals, the \emph{steps} $v = (v_1,\ldots,v_d)$ are a $d$-tuple of elements of $G$, and the \emph{symmetry group} $H$ is a finite subgroup of $G$.
\begin{itemize}
\item We define the \emph{image} $\Image(\P) \subset G$ to be the set
$$ \Image(\P) := H + [-N_1,N_1]_\Z \cdot v_1 + \ldots + [-N_d,N_d]_\Z;$$
thus $\Image(\P)$ is the sum of a subgroup and the image of a GAP.
\item For any $t > 0$, we define the \emph{dilate} $\P_t$ of $\P$ to be the coset progression $\P_t = (t N, v, d, H)$ formed by dilating all the dimensions by $t$ but keeping the symmetry group fixed.
\item We say that $\P$ is \emph{proper} if all the elements $h + n_1 v_1 + \ldots + n_d v_d$ for $n_i \in [-N_i,N_i]_\Z$ and $h \in H$ are distinct.  More generally, we say that $\P$ is \emph{$t$-proper} for some $t > 0$ if $\P_t$ is proper, and \emph{infinitely proper} if it is $t$-proper for all $t$ (i.e. the elements $v_1,\ldots,v_d$ are algebraically independent modulo $H$).
\item We define the \emph{size} of $\P$ to be $\size(\P) := |\Image(\P)|$.   Observe that if $\P$ is proper if and only if
$\size(\P) = |H| \prod_{i=1}^d (2 \lfloor N_i \rfloor + 1)$.
\end{itemize}
\end{definition}

\begin{remark} Of course, GAPs correspond to the special case $H = \{0\}$; also, finite subgroups of $G$ are essentially coset progressions of rank $0$.  Note that $H$ itself may require a large number of generators, but that this number has no bearing on the rank of $P$.
\end{remark}

Coset progressions are essential tools in the study of sum sets on arbitrary groups.  We mention two key (and closely related) theorems from \cite{GR} in this regard:

\begin{theorem}[Chang's theorem in an arbitrary group]\label{bog}[Section 5]\cite{GR} Let $A \subset G$
 be a non-empty finite set such that $|2A| \leq K|A|$ for some $K \geq 2$.
 Then there exists a proper coset progression $\P$ in $G$ of rank $O( K^3 \log K )$ and
size $\size(\P) \geq \exp( - O( K^3 \log^2 K) ) |A|$ such that
$\Image(\P) \subseteq 2A-2A$.
\end{theorem}

\begin{theorem}[Freiman's theorem in an arbitrary group]\label{frei}\cite[Theorem 1.1]{GR} Let $A \subset G$
be a non-empty finite set such that $|2A| \leq K|A|$ for some $K \geq 2$.
 Then there exists a coset progression $\P$ in $G$ of rank $O( K^4 \log K )$ and size at most $\exp( O( K^4 \log^2 K ) |A|)$ such that $A$ is contained in a translate of $\Image(\P)$.
\end{theorem}

One should compare these results to Theorems \ref{has-proper} and \ref{in-proper} respectively.
Alternate proofs of these results (with slightly weaker constants) can also be found in \cite[Theorem 5.48]{TVbook} and \cite[Theorem 5.44]{TVbook} respectively.

In Section \ref{torsion-sec} we establish the following generalization of Theorem \ref{prog-john}:

\begin{theorem}[John's theorem for coset progressions]\label{prog-john-coset} Let $\P$ be a symmetric coset progression of rank $d \geq 0$, and let $t \geq 1$.  Then there exists a $t$-proper symmetric coset progression $\Q$ of rank at most $d$,
such that we have the inclusions
$$ \Image(\P_{t'}) \subseteq \Image( \Q_{O(d)^{3d/2} t t'} )$$
for all $t' > 0$ and
$$ \Image(\Q_{tt'}) \subseteq \Image(\P_{t'})$$
for all $t' \geq 1$.
In particular, we have
$$ \Image(\Q) \subseteq \Image(\P) \subseteq \Image( \Q_{O(d)^{3d/2} t} ).$$
Furthermore we have the size bound
$$ t^{-d} 2^{-d^2-O(d \log d)} \size(\P) \leq \size(\Q) \leq \size(\P).$$
Finally, the symmetry group of $\Q$ contains that of $\P$, and if
$\P$ is not $1/2$-proper, then $\Q$ can be chosen to have rank at most $d-1$.
\end{theorem}

Note that when $\P$ lies in a torsion-free group, the symmetry group must be trivial, and so Theorem \ref{prog-john-coset} contains Theorem \ref{prog-john} as a special case.

By repeating the proof of Corollary \ref{prog-john} we obtain

\begin{corollary}[John's theorem for coset progressions, again]\label{prog-john-coset2}
Let $\P$ be a symmetric coset progression of rank $d \geq 0$, and let $t \geq 1$.  Then there exists a $t$-proper symmetric coset progression $\Q$ of rank at most $d$,
such that we have the inclusions
$$ \Image(\P_{t'}) \subseteq \Image( \Q_{t'} )$$
for all $t' > 0$ and
$$ \Image(\Q_{tt'}) \subseteq \Image(\P_{O(d)^{3d/2} t t'})$$
for all $t' \geq 1$.  In particular, we have
$$ \Image(\P) \subseteq \Image(\Q) \subseteq \Image( \P_{O(d)^{3d/2} t} ).$$
We also have the size bound
$$ \size(\P) \leq \size(\Q) \leq t^d O(d)^{3d^2/2} \size(\P).$$
Finally, if $\P_{d^{3d/2} t}$ is not proper, then $\Q$ can be chosen to have rank at most $d-1$.
\end{corollary}

\subsection{Sumsets in groups}

Now we consider the question of establishing John-type theorems for iterated sum sets $lA$ for $l$ large.
We will be interested here in the ``additively structured'' case when there is plenty of additive relations
between elements of $A$, either because $A$ is contained in a structured set such as a progression,
or because the iterated sumsets $lA$ are fairly small.  For instance, we have the following recent
result of Szemer\'edi and Vu.

\begin{theorem}[Sumsets in integers]\label{szv}\cite{sv} For every
integer $d \geq 1$ there exists $C, \eps > 0$ such that the following
statement holds: whenever $N \geq 1$, $l \geq 1$ and $A \subseteq [1,N]_\Z$
are such that $l^d |A| \geq CN$, then there exists a proper symmetric
GAP $\Q$ of rank $1 \leq d' \leq d$ and size at least $\eps l^{d'} |A|$ such that $lA$ contains a translate of $\Image(\Q)$.
\end{theorem}

For $d=1$, much more precise results in this direction are known: see \cite{sar}, \cite{freiman}, \cite{pom},
 \cite{lev}.  For further discussion of this result, including the sharpness
  of the various bounds, we refer the reader to  \cite{sv} and \cite[\S 12]{TVbook}.  An
alternate proof of this result can be found in \cite[Theorem 12.4]{TVbook}.
 For variants of this theorem when $l=2,3,4$ is small, see \cite{bourgain},
 \cite{fhr}, \cite{green-structure}, and Theorem \ref{bog} above.

Using the above machinery, together with the Freiman-type theorems of Green and Ruzsa, we can now generalize this statement to arbitrary additive groups with more explicit bounds and a slightly stronger statement (giving both an upper and lower containment for $lA$).

\begin{theorem}[John's theorem for iterated sumsets]\label{arbitrary} There exists a positive integer $C_1$ such
that the following statement holds: whenever $d \geq 1$, $l \geq
1$ and $A \subset G$ is a non-empty finite set such that $l^d |A|
\geq 2^{2^{C_1 d^2 2^{6d}}} |lA|$, then there exists a proper
symmetric
 coset progression $\Q$ of rank $0 \leq d' \leq d-1$ and size $\size(\Q) \geq 2^{-2^{C_1 d^2 2^{6d}}} l^{d'} |A|$ and $x, x' \in G$ such that
$$ x + \Image(\Q) \subseteq lA \subseteq x' + 2^{2^{C_1 d^2 2^{6d}}} \Image(\Q).$$
\end{theorem}

\begin{remark} The triple exponential dependence on $d$ is somewhat unsatisfactory; a single exponential would be more natural.  One exponential arises from the current best known bounds in Freiman-type theorems (see Theorem \ref{bog}).  Another arises from the need to ensure that a sumset of $\Image(\Q)$ can cover a sumset of $A$, which can temporarily exponentiate the rank of the progression. It may be that one of these exponential losses can be removed, or that two of them can be run ``in parallel'', reducing the total loss to a double exponential, but we will not attempt to do so here.
In the asymptotic limit $l \to \infty$, much more about the structure of $lA$ is known, for instance $|lA|$ is eventually a polynomial in $l$ \cite{nathanson-subgroup}, \cite{nr}.  The behavior of $lA$ for large $l$ is also closely connected to Theorems \ref{bog}, \ref{frei}; see \cite{GR} for further discussion.
\end{remark}

\begin{remark} When $A$ is symmetric, thus $A = -A$, and one has $0 \in A$ (or $l$ is even), then one can take $x=x'=0$ by exploiting the identity $2lA = lA - lA$ and the inclusions $lA \subset l'A$ when $l' \geq l$ and $0 \in A$.
\end{remark}

We prove this theorem in Section \ref{arbitrary-sec}.  Note that Theorem \ref{szv} emerges as a special case since in that case we have $|lA| \leq |l[1,N]_\Z| \leq lN$.
The $d=1$ case of this theorem already implies

\begin{corollary}[S\'ark\"ozy's theorem in an arbitrary group] There exists an absolute constant $C_2 > 0$ such that the following statement holds: whenever $l \geq 1$ and $A \subset G$ is an element of a finite additive group such
that $l|A| \geq C_2 |G|$, then $lA$ is a coset of the subgroup
generated by $A-A$.
\end{corollary}

\begin{proof}  Applying Theorem \ref{arbitrary} with $d=1$ we see that (for $C_2$ large enough)
we obtain a coset progression $\Q$ of rank $0$ with
 $\size(\P) \geq \eps l^{d'} |A|$ for some absolute constant $\eps > 0$,
 with $lA$ containing a translate of $\Image(\Q)$ and being contained in a translate
 of $C \Image(\Q)$ for some absolute constant $C$. Since $Q$ has rank $0$, $Q$ is a subgroup
and it is easy to check that $Q$ is generated by $A-A$. The claim
follows.
\end{proof}

\section{The discrete John's theorem}\label{djt-sec}

We now prove Theorem \ref{djt}, following the arguments in \cite[Section 3.5]{TVbook}.
We may assume that $\Gamma$ has full rank (i.e. its linear span is all of $\R^d$), since otherwise we can restrict to the linear span of $\Gamma$.

Applying John's theorem and a linear transformation we may assume that
\begin{equation}\label{bad}
B_d \subseteq B \subseteq \sqrt{d} \cdot B_d
\end{equation}
where $B_d$ is the (closed) unit ball in $\R^d$.  We recall the standard formula
\begin{equation}\label{ball-mes}
\mes(B_d) = \frac{\Gamma(3/2)^d 2^d}{\Gamma(d/2+1)} = \Theta(d)^{-d/2}
\end{equation}
where $\mes$ denotes $d$-dimensional Lebesgue measure.  In particular
\begin{equation}\label{b-mes}
\mes(B) \leq \sqrt{d}^d \mes(B_d) \leq O(1)^d.
\end{equation}

By the theory of Mahler bases (specifically, see \cite{bilu} or
\cite[Corollary 3.35]{TVbook}) we may find
linearly independent $v_1,\ldots,v_d$ which generate $\Gamma$, and such that
\begin{equation}\label{rmes}
\mes(\R^d/\Gamma) = |v_1 \wedge \ldots \wedge v_d| \geq O(d)^{-3d/2} |v_1| \ldots |v_d|.
\end{equation}
Indeed, the argument given in \cite{TVbook} allows one to write the $O(d)^{-3d/2}$ factor more explicitly as $\frac{\mes(B_d)}{2 d!}$.

A volume packing argument (see \cite[Lemma 3.26]{TVbook}) gives
\begin{equation}\label{bgam}
|B \cap \Gamma| \leq \frac{3^d d! \mes(B)}{2^d \mes(\R^d/\Gamma)} \leq \frac{O(d)^{5d/2}}{|v_1| \ldots |v_d|}
\end{equation}
where we have used \eqref{b-mes}.

Let us now set $\P$ to be the symmetric GAP $\P = (P, N, v, d)$, where $v := (v_1,\ldots,v_d)$, $N := (N_1,\ldots,N_d)$, and $N_i := \frac{1}{d |v_d|}$.  Since $v_1,\ldots,v_d$ are linearly independent, we see that $\P$ is infinitely proper. Observe that
$$ \size(\P) \geq \prod_{i=1}^d N_d = \frac{d^{-d}}{|v_1| \ldots |v_d|}$$
and so the first inequality in \eqref{dgam} follows from \eqref{bgam}.  Next, observe that if $x \in \Image(P_t)$ then
$$ |x| \leq tN_1 |v_1| + \ldots + tN_d |v_d| = t$$
and so $x \in (t \cdot B_d) \subseteq B$.  This gives the second containment in \eqref{inclusions-2}; setting $t=1$ we also obtain the second inequality in \eqref{dgam}.

Now let $x \in \Gamma$.  Since $v_1,\ldots,v_d$ generate $\Gamma$, we have $x = n_1 v_1 + \ldots + n_d v_d$ for some integers $n_1,\ldots,n_d$.  From Cramer's rule we have
$$ |n_i| = \frac{|x \wedge v_1 \wedge \ldots \wedge v_{i-1} \wedge v_{i+1} \wedge \ldots \wedge v_d|}{|v_1 \wedge \ldots \wedge v_d|}
\leq |x| \frac{|v_1| \ldots |v_d|}{|v_i| |v_1 \wedge \ldots \wedge v_d|}.$$
Applying \eqref{rmes} and the definition of $N_i$ we conclude that
$$ |n_i| \leq O(d)^{3d/2} |x| N_i$$
and so
$$ x \in \Image( \P_{O(d)^{3d/2} |x|} ).$$
Using this and \eqref{bad} we obtain the first inclusion in \eqref{inclusions-2}.  Since \eqref{inclusions-2} clearly implies \eqref{inclusions}, we obtain Theorem \ref{djt}.
\endprf

\section{Convex progressions}

It will be convenient to generalize from arithmetic progressions to the more geometric concept of
a \emph{convex progression}, as these are more stable under operations such as restriction or projection to
a subspace.  The idea of working with convex progressions was inspired by \cite{green}.  In view of our eventual generalization to the torsion case, we shall also allow the inclusion of a finite symmetry group, leading to \emph{convex coset progressions}

\begin{definition}[Convex progressions]
A \emph{symmetric convex coset progression} in an additive group $G$ is a quintuplet
$\P = (B, \Gamma, d, \phi, H)$, where the \emph{rank} $\rank(\P) := d \geq 0$ is an integer, $B$ is a convex body in $\R^d$, $\Gamma$ is a lattice in $\R^d$, $\phi: \Gamma \to G$ is a homomorphism.  We define the \emph{image} of $\P$ as
$$ \Image(\P) := \phi(B \cap \Gamma)+H \subset G$$
and the \emph{size} of $\P$ as $\size(\P) := |\Image(\P)|$.  If $t > 0$, we define the dilates $\P_t := (t \cdot B, \Gamma, d, \phi,H)$.  We say that $\P$ is \emph{proper} if the map $(x,h) \mapsto \phi(x)+h$ is is injective on $(B \cap \Gamma) \times H$, \emph{$t$-proper} if $\P_t$ is proper, and \emph{infinitely proper} if it is $t$-proper for all $t$.  If $H = \{0\}$ we refer to a symmetric convex coset progression as simply a \emph{symmetric convex progression}.
\end{definition}

Observe that every symmetric GAP $\P = (P, N, v, d)$ is also a symmetric convex progression of the same rank and size, setting\footnote{It is here that we are using our choice that convex bodies are compact rather than bounded open.  One could work with bounded open bodies by redefining GAPs using $(-N_i,N_i)_\Z$ instead of $[-N_i,N_i]_\Z$, but this creates some (minor) technical problems regarding dilations.} $B := \prod_{i=1}^d [-N_i,N_i]_\R$, $\Gamma := \Z^d$ and $\phi(n_1,\ldots,n_d) := n_1 v_1 + \ldots + n_d v_d$ for all $(n_1,\ldots,n_d) \in \Z^d$, and that the notions of dilation and properness are consistent.  Thus we can view convex progressions as generalizations of GAPs.  Similarly convex coset progressions are generalizations of coset progressions.

We observe the useful sumset embedding
$$\Image( \P_t ) +\Image( \P_{t'} ) \subseteq \Image(\P_{t+t'})$$
and hence
$$ k \Image(\P_t) \subseteq \Image(\P_{kt})$$
for $t,t' > 0$ and $k=1,2,\ldots$, and all symmetric convex coset progressions $\P$.
We shall use these embeddings frequently in the sequel without further comment.

\begin{remark} One way to view these embeddings is that they induce a translation-invariant pseudo-metric $\dist_\P: G \times G \to [0,+\infty]$ on $G$, defined by setting
$$ \dist_\P(x,y) := \inf \{ t > 0: x-y \in \Image(\P_t) \}$$
with $\dist_\P(x,y) = +\infty$ if no such $t$ exists (i.e. if $x-y$ is not in the group generated by $H$ and the $v_1,\ldots,v_d$).  We will however not use this metric structure explicitly.
\end{remark}

We now give two basic facts about these progressions.

\begin{lemma}[Doubling lemma]\label{cover} If $\P$ is a convex coset progression of rank $d$ and $t \geq 1$, then $\Image(\P_t)$ can be covered by at most $(4t+1)^d$ translates of $\Image(\P)$.  In particular
$$ \size(\P) \leq \size(\P_t) \leq (4t+1)^d \size(\P).$$
\end{lemma}

\begin{proof} It suffices to establish the first claim.  By applying the quotient map $G \to G/H$ we can reduce to the case $H=\{0\}$.  In this case the claim follows from \cite[Lemma 3.21]{TVbook}.  Here we use compact convex bodies rather than open bounded ones, but one can pass from one to the other by dilating by an epsilon and then sending that epsilon to zero; we omit the details.
\end{proof}

\begin{remark} In many cases one can lower the $(4t+1)^d$ factor to be closer to $t^d$, for instance if $t$ is an integer and $\P$ is a GAP with integer dimensions, but this type of improvement will not significantly improve our results here.
\end{remark}

\begin{lemma}[Covering lemma]\label{ruzsa} Let $\P, \Q$ be convex coset progressions of rank $d$ such that $\Image(\Q) \subseteq \Image(\P_t)$ for some $t > 0$.  Then for any $t' > 0$, $\Image(\P_{t'})$ can be covered by at most
$O(t+t'+1)^d \size(\P)/\size(\Q)$ translates of $\Image(\Q)$.
\end{lemma}

\begin{proof} Applying Ruzsa's covering lemma (see e.g. \cite[Lemma 2.14]{TVbook}) we conclude that
$\Image(\P_{t'})$ can be covered by $| \Image(\P_{t'})+\Image(\Q) |/|\Image(\Q)|$ translates of $\Image(\Q)-\Image(\Q)$.
But  by Lemma \ref{cover}
\begin{align*}
| \Image(\P_{t'})+\Image(\Q) | &\leq | \Image(\P_{t'})+\Image(\P_t) | \\
&\leq |\Image(\P_{t'+t})| \\
&=  \size(\P_{t'+t}) \\
&\leq O( t'+t+1)^d \size(\P)
\end{align*}
while $\Image(\Q)-\Image(\Q) \subseteq \Image(\Q_2)$ can be covered by $O(1)^d$ translates of $\Image(\Q)$.  The claim follows.
\end{proof}

\section{John's theorem in torsion-free groups}\label{pjt-sec}

The purpose of this section is to prove Theorem \ref{prog-john}.  The arguments here will be superceded by those in Section \ref{torsion-sec}, but we present these arguments first in the simpler torsion-free setting for expository purposes.

The key ingredient in the proof is the following rank reduction dichotomy in the convex progression setting.

\begin{lemma}[Lack of properness implies rank reduction]\label{lprop}  Let $\P$ be a symmetric convex progression of rank $d$ in a torsion-free group $G$, and suppose that $\P_{\frac{1}{2}}$ is not proper.  Then there exists a symmetric convex progression $\Q$ of rank $d-1$ with the inclusions
\begin{equation}\label{pt-include-1}
\Image(\P_t) \subseteq \Image(\Q_{2t})
\end{equation}
for all $t > 0$ and
\begin{equation}\label{pt-include-2}
\Image(\Q_t) \subseteq \Image(\P_t)
\end{equation}
for all $t \geq 1$.
\end{lemma}

\begin{proof} Write $\P = (B, \Gamma, d, \phi, \{0\})$.
 Since $\P_{1/2}$ is not proper, there exists distinct $x, x' \in (\frac{1}{2} \cdot B) \cap \Gamma$ such that $\phi(x) = \phi(x')$. Setting $y := x-x'$, we thus have $y \in B \cap \Gamma \backslash \{0\}$ and $\phi(y) = 0$.
 We may factorize $y = ny'$, where $n \geq 1$ is an integer and $y' \in \Gamma \backslash \{0\}$ is irreducible (thus $y' \neq m \cdot \Gamma$ for any integer $m>1$).  Since $G$ is torsion free, we thus have $y' \in B \cap \Gamma \backslash \{0\}$ and $\phi(y') = 0$.

Since $y'$ is irreducible, we can split $\Gamma = (\Z \cdot y') + \Gamma'$ where $\Gamma'$ is a lattice of rank at most $d-1$ (see e.g. \cite[Corollary 3.5]{TVbook}).  Applying an invertible linear transformation to $\Gamma$ and $B$ (and the inverse linear transformation to $\phi$) we can normalize $y' = e_d$ and $\Gamma' \subset \R^{d-1}$, where $e_1,\ldots,e_d$ is the standard basis of $\R^d$ and $\R^{d-1}$ is the span of $e_1,\ldots,e_{d-1}$.  Let $\pi: \R^d \to \R^{d-1}$ be standard projection, then we see that $\pi$ maps $\Gamma$ to $\Gamma' = \Gamma \cap \R^{d-1}$.  Since $\phi(e_d) = 0$, we see that we can factor $\phi = \phi' \circ \pi$ for some homomorphism $\phi': \Gamma' \to G$ (in fact $\phi'$ is just the restriction of $\phi$ to $\Gamma'$).  We now set
$\Q$ to be the symmetric convex progression
$$ \Q := (\frac{1}{2} \cdot \pi(B), \Gamma', d-1,\phi',\{0\}).$$
Since
$$ \phi( (t \cdot B) \cap \Gamma ) \subseteq \phi'( (t \cdot \pi(B)) \cap \Gamma' )
= \phi'( (2t \cdot \frac{1}{2} \cdot \pi(B)) \cap \Gamma' )$$
we obtain \eqref{pt-include-1}.  The only remaining thing to prove is \eqref{pt-include-2}, or in other words that
$$ \phi'( (\frac{t}{2} \cdot \pi(B)) \cap \Gamma' ) \subseteq \phi( (t \cdot B) \cap \Gamma ).$$
Since $\phi = \phi' \circ \pi$, it suffices to show that
$$ (\frac{t}{2} \cdot \pi(B)) \cap \Gamma' \subseteq \pi( (t \cdot B) \cap \Gamma ).$$
But if $z \in (\frac{t}{2} \cdot \pi(B)) \cap \Gamma'$, then $2z \in \pi(t \cdot B) \cap \Gamma'$, and so there exists $w \in \pi^{-1}(2z)$ such that $w \in t \cdot B$.  But since $e_d = y' \in B$, and $B$ is symmetric, we also have $t e_d, -t e_d \in t \cdot B$.  We conclude from convexity that $\frac{w}{2} + t' e_d \in B$ for all $-t/2 \leq t' \leq t/2$.  Since $\pi(w/2) = z \in \Gamma'$, $\Gamma = (\Z \cdot e_d) + \Gamma'$, and $t \geq 1$, we conclude that there exists
$-t/2 \leq t' \leq t/2$ such that $\frac{w}{2} + t' e_d \in (t \cdot B) \cap \Gamma$.  Since
$\pi(\frac{w}{2} + t' e_d) = z$, we obtain the desired inclusion.
\end{proof}

We can iterate the above lemma to obtain

\begin{corollary}[John's theorem for convex progressions]\label{conv-john}  Let $\P$ be a symmetric convex progression of rank $d \geq 0$ in a torsion-free group $G$, and let $t > 1/2$.  Then there exists a $t$-proper symmetric convex progression $\Q$ of rank $r$ for some $0 \leq r \leq d$, such that we have the inclusions
$$
\Image(\P_{t'}) \subseteq \Image( \Q_{2^{d-r+1} tt'} )$$
for all $t' > 0$ and
$$
\Image(\Q_{2tt'}) \subseteq \Image(\P_{t'})$$
for all $t' \geq 1$.  Furthermore, if $\P$ is not $1/2$-proper, then $\Q$ can be chosen to have rank $r \leq d-1$.
\end{corollary}

\begin{proof} It suffices to prove the claim for $t=1/2$, as the general case then follows by replacing $\Q$ by $\Q_{1/(2t)}$.  We now perform the following algorithm.
\begin{itemize}
\item Step 0.  Initialize $d=r$ and $\Q = \P$.
\item Step 1.  If $\Q_{1/2}$ is already proper then STOP.
\item Step 2.  Otherwise, use Lemma \ref{lprop} to replace $\Q$ by a symmetric convex progression $\Q'$ of rank $r-1$ such that $\Image(\Q') \subseteq \Image(\Q)$ and $\Image(\Q_t) \subseteq \Image(\Q'_{2t})$ for all $t > 0$.
\item Step 3.  Replace $r$ by $r-1$ and return to Step 1.
\end{itemize}
This algorithm terminates with $r \geq 0$ since progressions of rank $0$ are trivially proper.  The claims are then easily verified.
\end{proof}

\begin{proof}[Proof of Theorem \ref{prog-john}]
It suffices to verify the claim for $t = 1$, since the claim for larger $t$ then follows by
replacing $\Q$ by $\Q_{1/t}$ (and using Lemma \ref{cover} to recover the size bound for $\Q$).

Since $\P$ is a symmetric GAP, it is also a symmetric convex progression.  We can thus invoke Corollary \ref{conv-john} and find a convex $1/2$-proper progression $\P' = (B', \Gamma', d',\phi',\{0\})$ of rank $d' \leq d$ with the inclusions
$$ \Image(\P_t) \subseteq \Image(\P'_{2^d t}) \hbox{ for all } t > 0$$
and
$$\Image(\P'_t) \subseteq \Image(\P_t) \hbox{ for all } t \geq 1.$$
Also, if $\P$ is not $1/2$-proper, then $d' \leq d-1$.

By Lemma \ref{cover} we conclude
$$ \size(\P) \leq \size(\P'_{2^d}) \leq 2^{d^2+O(d)} \size(\P'_{1/2}).$$
Now, by Theorem \ref{djt} (applied to $\frac{1}{2} \cdot B'$ and $\Gamma'$)
we can find a proper symmetric GAP $\tilde \Q$ in $\R^{d'}$ of rank at most $d$ such that
$$  (O(d)^{-3d/2}t \cdot B') \cap \Gamma' \subseteq \Image(\tilde \Q_t) \subseteq (\frac{t}{2} \cdot B') \cap \Gamma'$$
for all $t > 0$, and also
$$ \size(\tilde \Q) \geq O(d)^{-7d/2} |\frac{1}{2} B' \cap \Gamma'|  O(d)^{-7d/2} \size(\P'_{1/2})$$
and thus
$$ \size(\tilde Q) \geq 2^{-d^2-O(d \log d)} \size(\P).$$
We can push forward the steps in $\tilde \Q$ by $\phi'$ to create a symmetric GAP $\Q$ in $G$ of rank at most $d' \leq d$.  Since $\tilde \Q$ is proper and $\phi'$ is injective on $(\frac{1}{2} \cdot B') \cap \Gamma$
we see that $\Q$ is also proper and has the same size as $\tilde Q$.  Now for any $t \geq 1$ we have
$$ \Image(\tilde \Q_t) \subseteq (\frac{t}{2} \cdot B') \cap \Gamma' \subseteq (t \cdot B') \cap \Gamma'$$
and thus applying $\phi'$ we have
$$ \Image(\Q_t) \subseteq \Image(\P'_t) \subseteq \Image(\P_t)$$
as desired.
Similarly, for any $t > 0$ we have
$$ (2^d t \cdot B') \cap \Gamma' \subseteq \Image(\tilde \Q_{O(d)^{3d/2} t})$$
and thus on applying $\phi'$ we have
$$ \Image(\P_t) \subseteq \Image(\P'_{2^d t}) \subseteq \Image(\Q_{O(d)^{3d/2} t})$$
as desired.  The claims of Theorem \ref{prog-john} have now all been established.
\end{proof}

\begin{remark}
Note that the bounds in Corollary \ref{conv-john} are slightly better than those for Theorem \ref{prog-john} (the dilations are of the order of $2^{O(d)}$ rather than $d^{O(d)}$, although the size bounds are comparable).  This suggests that in applications it may be more efficient to work with convex progressions instead of GAPs whenever possible.  See also \cite{green} for further discussion.
\end{remark}

\section{The torsion case}\label{torsion-sec}

We now extend the arguments of the previous section to the torsion case to prove Theorem \ref{prog-john-coset}.  The key ingredient is the following torsion variant of Lemma \ref{lprop}.

\begin{lemma}[Lack of properness implies rank reduction]\label{lprop2}  Let $\P$ be a symmetric coset convex progression of rank $d$ in a group $G$ (not necessarily torsion-free), and suppose that $\P_{\frac{1}{2}}$ is not proper.  Then there exists a symmetric convex coset progression $\Q$ of rank $d-1$ with the inclusions
\begin{equation}\label{pt-include-1a}
\Image(\P_t) \subseteq \Image(\Q_{2t})
\end{equation}
for all $t > 0$ and
\begin{equation}\label{pt-include-2a}
\Image(\Q_t) \subseteq \Image(\P_t)
\end{equation}
for all $t \geq 1$.  Furthermore, the symmetry group of $\Q$ contains the symmetry group of $\P$.
\end{lemma}

\begin{proof}  We follow closely the proof of Lemma \ref{lprop}.
Write $\P = (B, \Gamma, d, \phi,H)$.
Since $\P_{1/2}$ is not proper, there exists $x, x' \in (\frac{1}{2} \cdot B) \cap \Gamma$ and $h, h' \in H$ with
$(x,h) \neq (x',h')$ such that $\phi(x)+h = \phi(x')+h'$.  Note that this forces $x \neq x'$. Setting $y := x-x'$, we thus have $y \in B \cap \Gamma \backslash \{0\}$ and $\phi(y) \in H$.

As before we can find $y \in B \cap \Gamma \backslash \{0\}$ with
$\phi(y) \in H$, and can write $y=ny'$ where $y' \in B \cap \Gamma \backslash \{0\}$ is irreducible.
Unfortunately, we can no longer conclude that $\phi(y') \in H$, only that $n \phi(y') \in H$.  However, we can
shrink $n$ (and thus $y$) and assume that $n$ is the minimal positive integer such that $n \phi(y') \in H$; note
that we still have $ny \in B \cap \Gamma \backslash \{0\}$.  Let $H'$ be the group generated by $H$ and $\phi(y')$; this is then finite (indeed $|H'|=n|H|$).  In particular, $\phi(y')$ has finite order.

As before we can normalize $y'=e_d$ and split $\Gamma = (\Z \cdot e_d) + \Gamma'$ where $\Gamma' \subset \R^{d-1}$, and let $\pi: \R^d \to \R^{d-1}$ be the usual projection, thus $\pi$ maps $\Gamma$ to $\Gamma'$.  If we define $\phi': \Gamma' \to G$ to be the restriction of $\phi$ to $\phi'$, then it is no longer true that $\phi = \phi' \circ \pi$; instead, $\phi(x)$ and $\phi'(\pi(x))$ can differ by an element of the finite cyclic group
$\langle \phi(e_d) \rangle$.  We now set
$\Q$ to be the symmetric convex coset progression
$$ \Q := (\frac{1}{2} \cdot \pi(B), \Gamma', d-1,\phi', H').$$
Since
$$ \phi( (t \cdot B) \cap \Gamma ) + H \subseteq \phi'( (t \cdot \pi(B)) \cap \Gamma' ) + \langle \phi(e_d) \rangle + H
= \phi'( (2t \cdot \frac{1}{2} \cdot \pi(B)) \cap \Gamma' ) + H'$$
we obtain \eqref{pt-include-1a}.  The only remaining thing to prove is \eqref{pt-include-2}, or in other words that
$$ \phi'( (\frac{t}{2} \cdot \pi(B)) \cap \Gamma' ) + H' \subseteq \phi( (t \cdot B) \cap \Gamma ) + H.$$
Accordingly, let $z \in (\frac{t}{2} \cdot \pi(B)) \cap \Gamma'$ and $h' \in H'$; we need to show that
$$ \phi(z) + h' \in \phi( (t \cdot B) \cap \Gamma ) \mod H.$$
By definition of $H'$, we have $h' = l \phi(e_d) \mod H$ for some integer $l$.
As in the proof of Lemma \ref{lprop} we can find $w \in \pi^{-1}(2z)$ such that $w \in t \cdot B$.  But
$n e_d = y \in B$, and so by convexity as before we have $\frac{w}{2} + t' e_d \in B$ for all $-tn/2 \leq t' \leq tn/2$.
Since $\pi(w/2) = z$, we can express $w/2 = z + m e_d$ for some integer $m$, thus
$$ z + (m+t') e_d \in B \cap \Gamma \hbox{ for all } -tn/2 \leq t' \leq tn/2.$$
Since $t' \geq 1$, we may select $t'$ such that $m+t'-l$ is a multiple of $n$.  Since $n\phi(e_d) = 0 \mod H$, we conclude that
$$\phi(z) + h' = \phi(z + (m+t')e_d) - (m+t'-l) \phi(e_d) = \phi(z+(m+t')e_d) \mod H$$
and the claim follows.
\end{proof}

We can then iterate the proof of Corollary \ref{conv-john} more or less verbatim to obtain

\begin{corollary}[John's theorem for convex coset progressions]\label{conv-john2}  Let $\P$ be a symmetric convex coset progression of rank $d \geq 0$, and let $t > 1/2$.  Then there exists a $t$-proper symmetric coset convex progression $\Q$ of rank $r$ for some $0 \leq r \leq d$, such that we have the inclusions
$$
\Image(\P_{t'}) \subseteq \Image( \Q_{2^{d-r+1} tt'} )$$
for all $t' > 0$ and
$$
\Image(\Q_{2tt'}) \subseteq \Image(\P_{t'})$$
for all $t' \geq 1$.  Furthermore, the symmetry group of $\Q$ contains that of $\P$, and if $\P$ is
not $1/2$-proper, then $\Q$ can be chosen to have rank $r \leq d-1$.
\end{corollary}

Theorem \ref{prog-john-coset} then follows from Corollary \ref{conv-john2} in exactly the same way that Theorem \ref{prog-john} follows from Corollary \ref{conv-john} with only minor changes (e.g. replacing ``GAP'' with ``coset progression'' throughout, and tensor summing $\phi'$ with the identity
map on $H$).  We leave the details to the reader.

\section{Coalescence of progressions}

We now prove a variant of Theorem \ref{arbitrary} for coset progressions, which will in fact play a crucial role
in the proof of that theorem.

\begin{lemma}[Coalesence of coset progressions]\label{coalesce-2}  Let $\P$ be a symmetric coset progression of rank $d \geq 0$, and let $l \geq 1$ be an integer.  Then there exists a proper symmetric coset progression $\Q$ of rank $0 \leq d' \leq d$ such that
$\size(\Q) \geq d^{-3d^2/4} 2^{-C_0 d^2} l^{d'} \size(\P)$ for some absolute constant $C_0 > 0$
and
$$\Image(\Q) \subseteq l\Image(\P) \subseteq \lfloor (C_0 d)^{3d^2/4}\rfloor \Image(\Q).$$
In particular $\Image(\Q)$ generates the same group as $\Image(\P)$.
\end{lemma}

See \cite[Lemma 12.6]{TVbook} for a one-sided variant of this result in the torsion-free setting; earlier results in this direction are in \cite{sv}.  In fact our proof here is based on the proof of that lemma.

\begin{proof} We induct on $d$.  The case $d=0$ is trivial (setting $\Q :=\P$), so suppose that
$d \geq 1$ and the claim has already been proven for all smaller values of $d$.  We will fix $C_0$ to be a very large absolute constant to be chosen later.

By shrinking the dimensions we may assume that all dimensions of $\P$ are integers; in particular, $k\Image(\P) = \Image(\P_k)$ for all $k \geq 1$.  We may also assume without loss of generality that $l$ is a power of two and $l \geq (Cd)^{3d/2}$, where $C$ is a large absolute constant to be chosen later.

Let $l'$ be the greatest integer less than $l/(Cd)^{3d/2}$.
Suppose first that $\P$ is $l'$-proper.  Then we see that
$$ |l' \Image(\P)| = |\size(\P_{l'})| \geq O(Cd)^{-3d/2} l^d \size(\P)$$
and the claim follows by taking $\Q := \P_{l'}$. Thus we may assume that $\P$ is not $l'$-proper.  Thus there exists
$0 \leq k \leq \log_2 l'$ such that $\P$ is $2^k$-proper but not $2^{k+1}$-proper.  In particular, $\P_{2^{k+2}}$ is not $1/2$-proper.  Hence by Corollary \ref{prog-john2} (with $t=1$) we can find a proper symmetric coset progression $\P'$ of rank $0 \leq d'' \leq d-1$ such that
\begin{equation}\label{kball}
\Image(\P_{2^{k+2}}) \subseteq \Image(\P') \subseteq \Image( \P_{O(d)^{3d/2} 2^k} ).
\end{equation}
In particular we see that $\Image(\P')$ generates the same group as $\Image(\P)$, and also
\begin{equation}\label{deek}
\size(\P') \geq \size(\P_{2^k}) \geq 2^{dk - O(d)} \size(\P).
\end{equation}
We then apply the induction
hypothesis with $\P$ replaced by $\P'$ and $l$ replaced by $l'' = O(d)^{-3d/2} 2^{-k} l$ and conclude that there exists a proper symmetric coset progression $\Q$ of rank $0 \leq d' \leq d''$ such that
\begin{equation}\label{eek}
\Image(\Q) \subseteq l'' \Image(\P') \subseteq \lfloor (C_0 (d-1))^{3(d-1)^2/4}\rfloor \Image(\Q)
\end{equation}
and
$$ \size(\Q) \geq (d'')^{-3(d'')^2/4} 2^{-C_0 (d'')^2} (l'')^{d'} \size(\P').$$
Since $d'' \leq d-1$, we thus conclude (estimating $(d'')^{-3(d'')^2/4}$ from below by
$d^{-3(d-1)^2/4}$, and using \eqref{deek}) that (if $C_0$ is large enough)
$$ \size(\Q) \geq d^{-3d^2/4} 2^{-C_0 d^2} \size(\P).$$
Meanwhile, from \eqref{kball}, \eqref{eek} we have (if $l''$ was chosen correctly, and $C_0$ is large enough) that
\begin{align*}
\Image(\Q) &\subseteq l \Image(\P) \\
&\subseteq \lfloor (C_0 (d-1))^{3(d-1)^2/4}\rfloor  O(d)^{3d/2} \Image(\Q)\\
&\subseteq \lfloor (C_0 d)^{3d^2/4} \rfloor \Image(\Q)
\end{align*}
and the claim follows.
\end{proof}

\section{Sumsets in arbitrary groups}\label{arbitrary-sec}

We now prove Theorem \ref{arbitrary}.  We broadly follow the
strategy in \cite{sv}, showing that iterated sumsets first contain
large symmetric sets, then large coset progressions of high rank,
then large coset progressions of low rank.

We first translate $A$ so that $0 \in A$, so that the group
generated by $A-A$ is nothing more than $\langle A \rangle$, the
group generated by $A$. Also we see that the iterated sumsets $kA$
are nested in $k$. We shall also take $C_1$ to be a sufficiently
large absolute constant to be chosen later.

For any integer $0 \leq k \leq \log_2 l$, we have
$$ |2^k A| \leq |lA| \leq 2^{-2^{C_1 d^2 2^{6d} }} l^d |A|.$$
In particular this forces
\begin{equation}\label{l-large}
l \geq 2^{2^{C_1 d^2 2^{6d} }/d}.
\end{equation}
If we set $k'$ to be the first integer for which
\begin{equation}\label{2ka}
 |2^{k'+1} A| \leq 2^d |2^{k'} A|
\end{equation}
we thus see (from the pigeonhole principle) that $k'$ exists and is less than $\log_2 l$, and furthermore since
$$ |2^{k'} A| \geq 2^{dk'} |A|$$
we have
$$ 2^{dk'} \leq 2^{-2^{C_1 d^2 2^{6d} }} l^d$$
and thus in particular
\begin{equation}\label{kp-large}
 k' \leq \log_2 l - 2^{C_1 d^2 2^{6d}}/d.
\end{equation}
We will shortly encounter a need to replace $2^{k'} A$ by a symmetric set.  For this
we need the following lemma (a corrected version of \cite[Exercise 2.3.14]{TVbook}):

\begin{lemma}[Small doubling implies large symmetric sets] Let $A \subset G$ be finite with $|2A| \leq K|A|$ for some $K \geq 1$.  Then there exists $F \subseteq A$ and $x \in G$ such that $F = x-F$ and $|F| \geq |A|/K$.
\end{lemma}

\begin{proof} There are $|A|^2$ possible sums of the form $a+b$ with $a,b \in A$, which lie in a set $2A$ of cardinality at most $K|A|$.  By the pigeonhole principle, we can thus find $x \in 2A$ which can be written as such a sum in at least $|A|/K$ ways.  The claim then follows by setting $F := \{ a \in A: x-a \in A \}$.
\end{proof}

Applying this to $2^{k'} A$ we can find $F \subseteq 2^{k'} A$ and $x_0$ such that $F = x_0 - F$ and
$$ |F| \geq 2^{-d} |2^{k'} A| \geq 2^{dk'-O(d)} |A|.$$
In particular,
$$ |2F| \leq |2^{k'+1} A| \leq 2^{2d} |F|.$$
Applying Theorem \ref{bog} to $F$, we conclude that there exists a symmetric proper coset progression $\P$ of rank
$r = O( d 2^{6d} )$ and size
\begin{equation}\label{sizep}
\size(\P) \geq 2^{- O( d^2 2^{6d} )} |F| \geq 2^{- O( d^2 2^{6d} ) } |2^{k'} A|
\geq 2^{dk' - O( d^2 2^{6d} )} |A|
\end{equation}
such that
$$ \Image(\P) \subseteq 2F - 2F.$$
Since $F = x_0-F$, we have
$$ 2F-2F = 4F - 2x_0 \subseteq 2^{k'+2}A - 2x_0$$
and hence
$$ 2x_0 + \Image(\P) \subseteq 2^{k'+2} A.$$
In particular
$$ |2^{k'+2} A + \Image(\P)| \leq |2^{k'+3} A|.$$
But from \eqref{2ka} and Pl\"unnecke estimates (see \cite[Section 6.5]{TVbook}) we have
$$ |2^{k'+3} A| \leq 2^{O(d)} |2^{k'+2} A|.$$
Using the Ruzsa covering lemma as in Lemma \ref{cover}, we thus see that we can cover $2^{k'+2} A$ by
up to $2^{O(d^2 2^{6d})}$ translates of $\Image(\P)$.  Thus we may write
$$ 2^{k'+2} A \subseteq \bigcup_{i=1}^m 2x_0 + a_i + \Image(\P)$$
for some $a_1,\ldots,a_m$ with
\begin{equation}\label{msize}
m = O( 2^{O(d^2 2^{6d})} ),
\end{equation}
 which we may take to lie in
$$ a_i \in 2^{k'+2} A - (2x_0 + \Image(\P)) \in 2^{k'+2} A - 2^{k'+2} A.$$
Thus we can write $a_i = b_i - c_i$ for some $b_i, c_i \in 2^{k'+2} A$.

If $\P = (N, v, r, H)$, we let $\P'$ be the larger coset progression
$$ \P' := (N \oplus (1,\ldots,1), v \oplus (a_1,\ldots,a_m), r+m, H)$$
formed by adjoining $m$ new steps $a_1,\ldots,a_m$ with dimensions one each.  Then we have
\begin{equation}\label{psize}
\size(\P') \geq \size(\P) \geq 2^{dk' - O( d^2 2^{6d} )} |A|
\end{equation}
and
\begin{equation}\label{prank}
 \rank(\P') = r+m \leq 2^{O(d^2 2^{6d})}.
 \end{equation}
Also, by construction we have
\begin{equation}\label{2ka2}
2^{k'+2} A \subseteq \bigcup_{i=1}^m 2x_0 + a_i + \Image(\P) \subseteq \Image(\P').
\end{equation}
Furthermore, we observe that
$$ \Image(\P') \subseteq \Image(\P) + [-1,1]_\Z \cdot (b_1-c_1) + \ldots + [-1,1]_\Z \cdot (b_d - c_d).$$
Since $b_i,c_i \in 2^{k'+2} A$, we have $[-1,1]_\Z \cdot (b_i-c_i) \subseteq -b_i-c_i + 2^{k'+3} A$
and thus
$$ x_1 + \Image(\P') \subseteq (m+1) 2^{k'+3} A$$
where $x_1 := 2x_0 + \sum_{i=1}^m (b_i + c_i)$.

Let $l'$ be the largest integer such that $l'(m+1) 2^{k'+3} \leq l$; note from \eqref{l-large}, \eqref{kp-large} that $l'$ is at least $1$ (if $C_1$ is large enough).  Applying Lemma \ref{coalesce-2} we can thus find a symmetric coset progression $\Q$ of rank $0 \leq d' \leq r+m$ such that
\begin{equation}\label{drank}
\size(\Q) \geq (r+m)^{-O(r+m)^2} (l')^{d'} \size(\P') \geq 2^{-2^{O( d^2 2^{6d} )}} l^{d'} |A|
\end{equation}
(using \eqref{msize}, \eqref{psize}, \eqref{prank})
and
\begin{equation}\label{im}
 \Image(\Q) \subseteq l' \Image(\P') \subseteq O( (r+m)^{O((r+m)^2)} ) \Image(\Q) \subset
2^{2^{O(d^2 2^{6d})}} \Image(\Q).
\end{equation}
On the other hand, we have
$$ l' \Image(\P') \subseteq -l'x_1 + l' (m+1) 2^{k'+3} A \subseteq -l' x_1 + lA$$
and so $lA$ contains a translate of $\Image(Q)$.  On the other hand, from \eqref{2ka2}, \eqref{im} we have
$$ l' 2^{k'+2} A \subseteq 2^{2^{O(d^2 2^{6d})}} \Image(\Q)$$
and hence (by definition of $l'$)
$$ l A \subseteq O( m 2^{2^{O(d^2 2^{6d})}}) \Image(\Q) \subseteq 2^{2^{O(d^2 2^{6d})}} \Image(\Q).$$

We now have the required size and containment bounds on $\Q$.  The only problem is that we have
only a very poor bound on $d'$.  We can improve it by noting that
$$ \size(\Q) \leq |lA| \leq 2^{2^{-C_1 d^2 2^{6d}}} l^d |A|$$
and hence by \eqref{drank} we have $d' \leq d-1$ if $C_1$ is
sufficiently large (compared to $d$ and the hidden constant in the
$O$ of \eqref{drank}). The proof of Theorem \ref{arbitrary} is now
complete.
\endprf

\end{document}